\documentclass[ 12pt]{article} 

\usepackage[english]{babel}
\usepackage{amsmath,amssymb,amsthm}
\usepackage{mathtools}
\usepackage{microtype}
\usepackage{hyperref} % keep hyperlinks; the production team will adjust
\usepackage{graphicx}
\usepackage{booktabs}
\usepackage{enumitem}
\usepackage{cleveref}

\numberwithin{equation}{section}   
\theoremstyle{plain}
\newtheorem{theorem}{Theorem}[section]    
\newtheorem{lemma}[theorem]{Lemma}

\theoremstyle{definition}

\theoremstyle{remark}

\newtheorem{notation}[theorem]{Notation}

\begin{document}

\title{Bounded Solutions of a Complex Differential Equation for the Riemann Hypothesis}
 
\author{W. Oukil\\
\small\text{Faculty of Mathematics.}\\
\small\text{University of Science and Technology Houari Boumediene.}\\
\small\text{ BP 32 EL ALIA 16111 Bab Ezzouar, Algiers, Algeria.}}

\date{\today}

\maketitle

%%%%%%%%%%%%%%%%%%%%%%%%%%%%%%%%%%%%%%%%%%%%%%%%%%%%%%%%%%

%%%%%%%%%%%%%%%%%%%%%%%%%%%%%%%%%%%%%%%%%%%%%%%%%%%%%%%%%%

\begin{abstract}
In this manuscript, we consider the Riemann zeta function $\zeta$, defined through the Abel summation formula. 
We present a simple analytical method based on a complex differential equation.

The aim is to propose a new analytical approach, relying on complex differential equations defined on the interval $[1,+\infty)$, 
in order to gain insight into the behavior of $\zeta(s)$ within the critical strip. 
We introduce a differential equation depending only on the complex parameter $s$, extracted from the analytical structure of $\zeta(s)$ for $s$ in the critical strip. 
This equation admits a unique continuous and bounded solution. 
The non-trivial zeros of the zeta function can thus be characterized through the boundedness of such a solution. 
Furthermore, we conjecture an asymmetry in the boundedness of these solutions with respect to the critical line, 
suggesting that if $\zeta(1-s)= 0$, then $\zeta(s) \neq 0$ for any $s$ in the critical strip except on the critical line. 
This observation does not contradict the Riemann functional equation but supports a formulation consistent with the Riemann Hypothesis, 
opening a simple yet potentially new direction for the analytical investigation of the zeta function and the localization of its non-trivial zeros.
\end{abstract}

\begin{keywords} 
Differential Equations, Riemann  zeta function,   Non trivial zeros,  Riemann functional equation, Riemann hypothesis.
\end{keywords}\\
\begin{MSC}
00A05,   34A30, 34E05  
\end{MSC}

%%%%%%%%%%%%%%%%%%%%%%%%%%%%%%%%%%%%%%%%%%%%%%%%%%%%%%%%%%

%%%%%%%%%%%%%%%%%%%%%%%%%%%%%%%%%%%%%%%%%%%%%%%%%%

%%%%%%%%%%%%%%%%%%%%%%%%%%%%%%%%%%%%%%%%%%%%%%%%%%%%%%%%%%

%%%%%%%%%%%%%%%%%%%%%%%%%%%%%%%%%%%%%%%%%%%%%%%%%%
\section{Introduction}  
Consider  the representation of the Riemann  zeta function $\zeta$  defined    by   Abel's   summation formula [\cite{Titchmarsh}, page 14  Equation 2.1.5]  as
\begin{equation}\label{Zetafunction}
\zeta(s):= -\dfrac{s}{1-s}-s\int_1^{+\infty}u^{-1-s}\left\{u\right\}du,\quad   s \in \mathbb{C},\quad s\neq1,\quad\Re(s)>0,
\end{equation}
where $\left\{u\right\}$ is the fractional part of the real $u$.  The core idea of this work lies in proposing a dynamical approach to the Riemann Hypothesis. The main novelty is to reformulate the Riemann Hypothesis as a question of stability and asymptotic behavior within the framework of ordinary differential equations. Specifically, the paper establishes that a complex number $s$ is a non-trivial zero of the Riemann zeta function if and only if a particular solution of the associated complex differential equation remains bounded. This equivalence translates the problem of zero distribution into a condition on the qualitative behavior (boundedness) of a dynamical system. 
 
\section{Zeta function as a complex differential equation} \label{SecMainTheorem}
\noindent
Denote by  $B\subset \mathbb{C}$  the critical strip, defined as
\[
{B}:=\Big\{s\in \mathbb{C}:\quad        \Re(s)\in(0, 1)\Big\}.
\]  
Let  $w\in {B}$.  Based on equation \eqref{Zetafunction}, the aim is to study  the differential equation of solutions the following functions
\[
t\mapsto \psi_{w}(z,t):=t^{w  }\Big[z+ \int_1^{t} u^{ -1-w}  \left\{u\right\} du\Big], \quad z\in\mathbb{C},\ t\geq1.
\]
The aim is to show that the point  $s\in B$ is a   zeros of the function $\zeta$ then it is  loaclized through the boundedness of the solution  $ t\mapsto \psi_w ( \frac{1}{1 - w}, t )$ on $[1,+\infty)$. In other words, $\sup_{t \geq 1}  | \psi_w ( \frac{1}{1 - w}, t  )   | < +\infty$ if and only if $|\frac{  \zeta(w)}{w}| =0$. That gives the following equivalence: The Reimann hypothesis is true if for every $s\in B$ such  $\Re(s)\neq \frac{1}{2}$ the functions $\psi_s( \frac{1}{1 - s}, \cdot )$ and $\psi_{1 -  {s}}( \frac{1  }{ {s}}, \cdot )$ can not be  both bounded on $[1, +\infty)$. In this paper we derive the functions only on $[1,+\infty)/\mathbb{N}$. The following lemma defines the differential equation in question.
\begin{lemma}\label{lemmeIntegrationExpli}
For every  $w\in {B}$ and $z\in \mathbb{C}$ there exists a unique continuous solution $\psi_{{w}}(z, \cdot )$       of      the  following differential equation 
\begin{gather}
\notag     \dfrac{d}{dt} x = w   \, t^{-1}x+  \, t^{-1} \left\{t\right\},\\
\notag  \quad t\in [1,+\infty)/\mathbb{N},\quad x( 1) =z, \quad x:[1,+\infty)\to\mathbb{C}.
\end{gather}  
Further,  
\[
 \psi_{{w}}(z,t)=  t^{{w} }  \Big[z+ \int_1^{t} u^{ -1-w}  \left\{u\right\} du\Big],\quad \forall t\geq1.
\]  
\end{lemma}
\begin{proof}
Let $w\in {B}$ fixed. The function
\[
t\mapsto  t^{{w}} \int_1^t   u^{-1-w} \left\{u\right\} du,\ t\geq1,
\]
is   $C^\infty$ on $[1,+\infty)/\mathbb{N}$ and continuous on $[1,+\infty)$. The differential equation in question is   a non-homogeneous linear differential equation. For every $z\in \mathbb{C}$ the unique continuous solution $\psi_w(z, \cdot )$ such that $\psi_{{w}}(z,1)=z$ is given by 
\[
 \psi_{ {w}}(z,t)  =t^{{w} } \Big[z+ \int_1^{t} u^{-1-w}  \left\{u\right\} du\Big], \quad \forall t\geq1.
\]
\end{proof}  
\begin{notation}\label{NotationPsi}
For every $w\in {B}$ and $z\in\mathbb{C}$, we denote the function $\psi_{w}(z, \cdot )$ from $[1,+\infty)$ to $\mathbb{C}$ the unique continuous solutions of the  following differential equation 
\begin{gather}
\label{equdiff}      \dfrac{d}{dt} x =w    \, t^{-1}x+   \, t^{-1} \left\{t\right\},\\
\notag  \quad t\in [1,+\infty)/\mathbb{N},\quad x( 1) =z, \quad x:[1,+\infty)\to\mathbb{C}.
\end{gather}  
The point  $z$ is called the initial condition of the solution  $ \psi_{w}(z, \cdot )$. 
\end{notation}
\begin{lemma}\label{perlemma} Denote
\begin{equation*}
p(u):= \int_1^{u}\left(\dfrac{1}{2}-\left\{v\right\} \right)dv,\quad \forall u\geq1.
\end{equation*}
For all $u\geq1$, we have  $0\le p(u)\le  \frac{1}{8}$.
\end{lemma}
\begin{proof}
It is sufficient to prove that the function $u\mapsto p(u)$ is 1-periodic. We have 
\[
p(u+2) = \int_1^{2}\left(\dfrac{1}{2}-\left\{v\right\} \right)dv+ \int_2^{2+u}\left(\dfrac{1}{2}-\left\{v\right\} \right)dv,\quad \forall u\geq0.
\]
Since
\[
  \int_1^{2}\left(\dfrac{1}{2}-\left\{v\right\} \right)dv=\int_1^{ 2}\left(\dfrac{1}{2}-v+1 \right)dv=0,    
\]
we get
\[
p(u+2) =   \int_2^{2+u}\left(\dfrac{1}{2}-\left\{v\right\} \right)dv, \quad u\geq0.
\]
Using the fact that the function  $v\mapsto\dfrac{1}{2}-\left\{v\right\}$ is  1-periodic,  we obtain
\[
p(u+2) = \int_1^{1+u}\left(\dfrac{1}{2}-\left\{v\right\} \right)dv=p(1+u), \quad u\geq0.
\]
Then the function  $u\mapsto p(u)$ is   1-periodic and  for all $u\geq1$ we get $0=\min_{v\in[1,2]}p(v)\le p(u)\le \frac{1}{8}=\max_{v\in[1,2]}p(v)$
\end{proof} 
In the following lemma, we are interested in the bounded solutions on $[1,+\infty)$ .
\begin{lemma}\label{lemmeIntegrationExpliB}
Let $w\in {B}$ and $z\in \mathbb{C}$. Consider  the solution $\psi_{w}(z , \cdot )$       of       the   differential equation \eqref{equdiff} as given in the Notation \ref{NotationPsi}. Then  $\sup_{t\geq1}|\psi_{ w}(z ,t )|<+\infty$ if and only if $z=z_w$, where
\[
z_w=- \int_1^{+\infty}  u^{-1-w} \left\{u\right\} du.
\]  
Further, there exists $c_w>0$ such that
\[
\left|\psi_{ w}(z_w,t )+\dfrac{1}{2 w}\right|<c_w   \, t^{-1},\quad \forall t\geq1.
\]
\end{lemma}
\begin{proof}
Let $w\in {B}$ fixed. There is a unique bounded solution. All other solutions are oscillating  and diverge to infinity in norm.
Let $z\in \mathbb{C}$ and suppose that $\sup_{t\geq1}|\psi_{ w}(z,t)|<+\infty$. By the Lemma \ref{lemmeIntegrationExpli}, we have
\[
\sup_{t\geq1}|\psi_{ w}(z,t)|<+\infty\implies\sup_{t\geq1}\Bigg[  t^{\Re(w)}\left|   z+ \int_1^{+\infty}  u^{-1-w} \left\{u\right\} du\right|\Bigg]<+\infty,
\]  
Since $\Re(w)>0$, then
\[
 \left|   z+ \int_1^{+\infty}  u^{-1-w} \left\{u\right\} du\right|=0,
\]  
by consequence, $z=z_w$. Prove that $\sup_{t\geq1}|\psi_{ {w}}(z_{{w}},t )|<+\infty$. By definition of the point $z_w$ we have
\[
  z_{{w}}+ \int_1^t   u^{-1-w} \left\{u\right\} du=- \int_t^{+\infty}  u^{-1-w} \left\{u\right\}  du ,\quad \forall t>1.
\]
By the Lemma \ref{lemmeIntegrationExpli}, we obtain
\begin{equation}\label{psiparticuliere}
 \psi_{{w}}(z_{{w}},t)=-  t^{{w} }  \int_t^{+\infty}  u^{-1-w} \left\{u\right\}  du,\quad \forall t\geq1.
\end{equation}
Then $\sup_{t\geq1}|\psi_{ {w}}(z_{{w}},t )|\le \dfrac{1}{\Re(w)}$.  Now, prove the second item of the present Lemma. From   equation\eqref{psiparticuliere}, we have
\[
 \psi_{{w}}(z_{{w}},t)=-\dfrac{1}{2w}  +  t^{{w} }\int_t^{+\infty}  u^{-1-w} (\dfrac{1}{2}-\left\{u\right\})  du,\quad \forall t\geq1.
\]
Using the integration by parts formula, we obtain
\[
 \psi_{{w}}(z_{{w}},t)=-\dfrac{1}{2w} +  (1+w)  t^w\int_t^{+\infty}  u^{-2-w}\int_t^{u}\left(\dfrac{1}{2}-\left\{v\right\} \right)dv du. 
\]
By   Lemma \ref{perlemma}, for every $t\geq1$, we get
\[
\left| \psi_{{w}}(z_{{w}},t)+\dfrac{1}{2w }   \right|\le    \, t^{-1} \dfrac{|1+w  |}{1+\Re(w)}.
\]
\end{proof}
\section{Main Theorem}
The core equivalence between the zeros of the Riemann zeta function and the boundedness of the associated differential equations solutions is stated below.
\begin{theorem}\label{BoundedSolutionandZeroZeta}
We have the following equivalences:
\[
\forall\, w \in B: \quad 
\zeta(w) = 0 
\iff  
\sup_{t \geq 1} \left| \psi_{w}\!\left(\tfrac{1}{1-w}, t\right) \right| < +\infty,
\]
and
\[
\forall\, w \in B: \quad 
\zeta(w) = 0 
\iff  
\sup_{t \geq 1} \left|\, t \!\left( \psi_{w}\!\left(\tfrac{1}{1-w}, t\right) + \dfrac{1}{2w} \right) \!\right| < +\infty.
\]
\end{theorem}

\begin{proof}
From Equation~\eqref{Zetafunction}, for every $w \in B$ we have
\begin{equation}
\label{zetaproof}
\dfrac{\zeta(w)}{w}
  = -\dfrac{1}{1-w}
    - \int_1^{+\infty} u^{-1-w} \{u\}\, du.
\end{equation}
Consider the bounded solution $\psi_{w}(z_w, \cdot)$ of the differential equation~\eqref{equdiff}, as given in Lemma~\ref{lemmeIntegrationExpliB}. 
For this solution, we have
\[
z_w = -\int_1^{+\infty} u^{-1-w} \{u\}\, du.
\]
Consequently, Equation~\eqref{zetaproof} can be rewritten as
\[
\dfrac{\zeta(w)}{w} = -\dfrac{1}{1-w} + z_w.
\]
Hence, $\zeta(w) = 0$ if and only if $z_w = \tfrac{1}{1-w}$. 
By uniqueness of the continuous solution, we then obtain 
\[
\psi_{w}\!\left(\tfrac{1}{1-w}, \cdot\right) = \psi_{w}(z_w, \cdot)
\quad \text{for all } t \geq 1.
\]
\end{proof}
\section{Conjecture}

We conjecture that if $\Re(s) \neq \tfrac{1}{2}$, then the continuous solutions 
$\psi_{s}\!\left(\tfrac{1}{1-s}, \cdot \right)$ and 
$\psi_{1-  \bar{s}}\!\left(\tfrac{1}{ \bar{s}}, \cdot \right)$ 
of the differential equation~\eqref{equdiff} cannot both be bounded on $[1,+\infty)$. 
Numerical computations, not included in this manuscript, are consistent with this conjecture.

\medskip
A more general conjecture can be stated when the fractional part function in the non-homogeneous linear differential equation~\eqref{equdiff} 
is replaced by a real, bounded, locally integrable function 
$\eta : [1,+\infty) \to \mathbb{R}$ satisfying
\[
  \sup_{u \geq 1} 
  \Bigg| 
    \int_{1}^{u} \!\Big(\tfrac{1}{2} - \eta(v)\Big) \, dv 
  \Bigg| < \tfrac{1}{2}.
\] 
From the viewpoint of the qualitative theory of differential equations,  Equation~\eqref{equdiff} defines a family of linear non-autonomous equations parametrized by $\Re(w)$.  The parameter $\Re(w)$ induces an ordering on the family of solutions, and one may study the dependence on this parameter by differentiating the solutions with respect to $\Re(w)$ 
and analyzing the corresponding variational equation.

We stress, however, that although the conjecture extends to this broader class of inhomogeneities, this fact alone is insufficient to deduce that the zero  lie on the critical line. In the case of the fractional part function, the additional symmetry provided by the Riemann functional equation is essential for drawing conclusions about the location of the zeros. The theoretical framework developed for this notion of inhomogeneity applies, in particular, to the standard generalizations of the zeta function, namely the family of $L$-functions.

When $\Re(s) \neq \tfrac{1}{2}$ and the continuous solutions  
$\psi_{s}\!\left(\tfrac{1}{1-s}, \cdot \right)$ and  
$\psi_{1-  \bar{s}}\!\left(\tfrac{1}{  \bar{s}}, \cdot \right)$  
of the differential equation~\eqref{equdiff} cannot both be bounded on $[1,+\infty)$,  
Theorem~\ref{BoundedSolutionandZeroZeta} implies that the two quantities  $\zeta(s)$ and $\zeta(1- \bar{s})$ cannot both vanish at the same time.  
A direct consequence, obtained by applying the Riemann functional equation~\cite[Theorem~2.1, p.~13]{Titchmarsh},  
is that the non-trivial zeros of the Riemann zeta function must lie on the critical line.

\makeatother

%%%%%%%%%%%%%%%%%%%%%%%%%%%%%%%%%%%%%%%%%%%%%%%%%%%%%  

\end{document}